\documentclass[preprint,12pt]{elsarticle}

\usepackage{amssymb}
\usepackage{amsthm}

\usepackage[dvipsnames]{xcolor}
\usepackage{amsmath,amssymb,amsfonts}
\usepackage{amsthm}
\newtheorem{theorem}{Theorem}[section]

\usepackage{algorithm}
\usepackage{algpseudocode}
\usepackage{color,soul}
\usepackage{cleveref}
\usepackage{comment}

\usepackage[listings,skins,breakable]{tcolorbox}
	{%
	\end{tcolorbox}\endgroup}

\newtcbox{\popovnotesmall}{breakable,enhanced jigsaw,nobeforeafter,tcbox raise base,boxrule=0.4pt,top=0mm,bottom=0mm,
	right=0mm,left=4mm,arc=1pt,boxsep=2pt,before upper={\vphantom{dlg}},
	colframe=Fuchsia!75!black,coltext=black,colback=Fuchsia!20,
	overlay={\begin{tcbclipinterior}\fill[Fuchsia!75!black] (frame.south west)
			rectangle node[text=white,font=\sffamily\bfseries\tiny,rotate=90] {AAP} ([xshift=4mm]frame.north west);\end{tcbclipinterior}}}
\MakeRobust\popovnotesmall
\ifdefined
\fi

\def\*#1{\boldsymbol{\mathbf{#1}}}
\DeclareMathOperator{\Cov}{Cov}
\definecolor{highlighter}{cmyk}{0,0,0.5,0}

\journal{Journal of Computational Physics}

\begin{document}
	
	\begin{frontmatter}
		
		
		
		\title{Bayesian Recursive Update for Ensemble Kalman Filters}
		
		
		\author[ase]{Kristen Michaelson}
		\affiliation[ase]{organization={Dept. of Aerospace  Engineering \& Engineering Mechanics, The University of Texas at Austin},
			addressline={2617 Wichita St. C0600}, 
			city={Austin}, 
			postcode={78712},
			state={Texas},
			country={United States}}
		
		\author[oden]{Andrey A. Popov}
		\affiliation[oden]{organization={Oden Institute for Computational Engineering and Sciences},
			addressline={201 E 24th St. C0200}, 
			city={Austin},
			postcode={78712}, 
			state={Texas},
			country={United States}}
		
		\author[ase]{Renato Zanetti}

		\begin{abstract}
			Few real-world systems are amenable to truly Bayesian filtering; nonlinearities and non-Gaussian noises can wreak havoc on filters that rely on linearization and Gaussian uncertainty approximations. This article presents the Bayesian Recursive Update Filter (BRUF), a Kalman filter that uses a recursive approach to incorporate information from nonlinear measurements. The BRUF relaxes the measurement linearity assumption of the Extended Kalman Filter (EKF) by dividing the measurement update into a user-defined number of steps. The proposed technique is extended for ensemble filters in the Bayesian Recursive Update Ensemble Kalman Filter (BRUEnKF). The performance of both filters is demonstrated in numerical examples, and new filters are introduced which exploit the theoretical foundation of the BRUF in different ways. A comparison between the BRUEnKF and Gromov flow, a popular particle flow algorithm, is presented in detail. Finally, the BRUEnKF is shown to outperform the EnKF for a very high-dimensional system. 
		\end{abstract}
		
		\begin{graphicalabstract}
			\includegraphics[width=\textwidth]{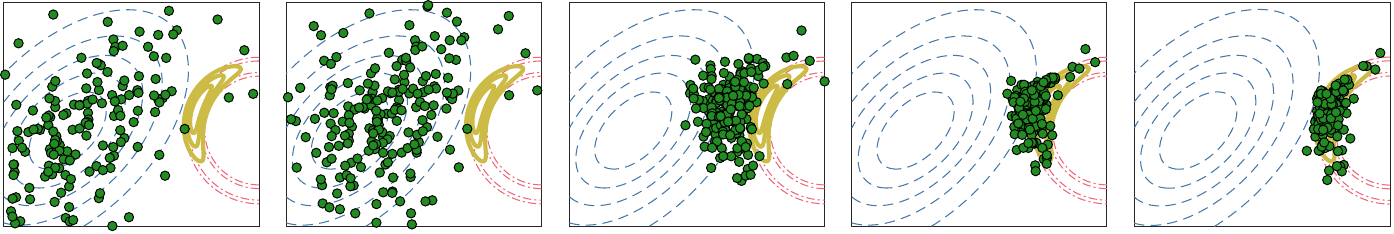}
		\end{graphicalabstract}
		
		\begin{highlights}
			\item The Bayesian Recursive Update Filter (BRUF) is presented and proven to match the result of the Kalman update for linear measurements.
			\item Two new forms of the BRUF are introduced: the variable-step BRUF (VS-BRUF), and the error controller BRUF (EC-BRUF). These new forms take advantage of a key insight into the linear convergence proof.
			\item The BRUF, the VS-BRUF, and the EC-BRUF are all adapted for use in an ensemble Kalman filter (EnKF). The performance of this new class of EnKFs is demonstrated on two high-dimensional problems.
		\end{highlights}
		
		
	\end{frontmatter}
	

	\section{Introduction}
	State estimation is the process of recursively estimating the state of a system given a dynamics model and a set of measurements. The state estimate is propagated forward in time using the dynamics model. As measurements arrive, the state estimate is updated to reflect the new information. If the noises in the dynamics and the measurements are Gaussian-distributed and the dynamics and measurement models are linear, then the Kalman filter provides the optimal state estimate in a mean-squared error sense \cite{Kalman_1960}. 
	
	If the dynamics and measurement models are nonlinear, the Extended Kalman Filter (EKF) approximates the optimal solution using linearizations of the dynamics and measurement models about the current state estimate \cite{bar2001estimation}. The linearization is used to propagate the estimated covariance; the state dynamics can still be propagated nonlinearly using numerical integration. However, there is not an obvious solution for adjusting the state estimate in a nonlinear way given new \textit{measurement} information. 
	The measurement update occurs at a single time step, and it must provide a solution that approximates the true Bayesian posterior distribution of the state. 
	High-precision nonlinear measurements pose a particular challenge to linear filters. The primary focus of this work is improving on the EKF approximation of the Bayesian posterior distribution for nonlinear measurements. 
	
	One approach to improving the EKF update for nonlinear measurements is the Iterated EKF (IEKF) \cite{gelb1974applied}. The IEKF begins by computing an EKF update. Then, the measurement model is linearized about the updated state. A new EKF update is computed using the new measurement linearization, and so on, until there is little change in the updated state between iterations. The covariance update is applied after the final iteration, using the final measurement linearization. There are no intermediate covariance updates. The IEKF state estimate has been shown to converge to the \textit{maximum a posteriori} (MAP), the state value of highest probability in the Bayesian posterior distribution \cite{bell1993iterated}. However, there is no convergence guarantee and a line search to insure descent is sometimes needed \cite{zanetti2011recursive}. 
	
	The recursive update filter (RUF) sought to smooth the EKF update by moving the state estimate through a series of attenuated Kalman updates \cite{zanetti2011recursive, zanetti2015adaptable}. The RUF also takes advantage of relinearization; however, instead of computing a series of state updates and only applying the last one, the RUF applies each state update successively. The Bayesian recursive update filter (BRUF), introduced in \cite{bruf}, is a similar update scheme that uses inflated measurement noise covariance instead of scaling the Kalman gain to attenuate the update steps. 
	
	In \cite{bruenkf}, the BRUF technique was extended for ensemble Kalman filters \cite{evensen2009data}. Ensemble filters and particle filters carry many state vectors in order to faithfully represent the probability distribution of the state. The Bayesian recursive update ensemble Kalman filter (BRUEnKF) measurement update is reminiscent of existing particle flow filters \cite{daum2010exact, daum2016gromov, ding2012implementation, crouse2019consideration}. 
	
	Similarities between ensemble filters and particle filters are easy to find. They both rely on a group of state estimates to represent an uncertainty distribution. Differences lie only in implementation. Ensemble filters traditionally employ an \textit{inflation} step before the measurement update \cite{popov2020explicit}. They are popular in the geoscience community, where solving systems with very large covariance matrices for very large state representations is intractable. The state estimate at any time step is simply the mean of the ensemble members, and necessary covariance values can be computed at measurement time using the statistical covariance of the Jacobian vector products of the difference between the ensemble and the mean. Particle filters, on the other hand, traditionally include scalar weight values for each state vector. They also often use a resampling step to redraw a new set of states after the measurement update \cite{ristic2003beyond}.  Particle \textit{flow} filters also include an optional resampling step, but they forego weights in favor of a companion EKF or unscented Kalman filter (UKF). 
	
	This article builds on two conference papers: the BRUF \cite{bruf} was first presented at at the 33rd AAS/AIAA Space Flight Mechanics Meeting, and the BRUEnKF \cite{bruenkf} was introduced at the 26th International Conference on Information Fusion. In this work, we provide new examples to demonstrate the performance of the BRUF and BRUEnKF updates, both for a single measurement update step and in filters. We also introduce two new versions of the BRUF and the BRUEnKF. These new versions take advantage of a key insight into the proofs of linear convergence for the BRUF. 
	
	The remainder of this paper is organized as follows: Section \ref{sec:bruf} presents the BRUF and a convergence proof for linear measurements; Section \ref{sec:range-observation-bruf} presents a 2D example of the BRUF update and introduces the variable-step BRUF (VS-BRUF) and the error-controller BRUF (EC-BRUF); Section \ref{sec:tracking-bruf} compares the performance of the proposed filters to the IEKF in an aerospace tracking scenario; Section \ref{sec:bruenkf} revisits the 2D example of Section \ref{sec:range-observation-bruf}, this time with the proposed ensemble filters. The ensemble filters are demonstrated in two examples: a Lorenz '96 system with a highly nonlinear measurement in Section \ref{sec:bruenkf-l96} and a geoscience example---the quasi-geostrophic equations---in Section \ref{sec:bruenkf-qg}. Finally, conclusions and future work are presented in Section \ref{sec:conclusion}. 
	
	\section{The Bayesian Recursive Update Filter}
	\label{sec:bruf}
	The BRUF 
	relaxes the measurement linearity assumption of the EKF by recursively applying EKF updates with inflated measurement noise. Consider a nonlinear measurement model with additive unbiased Gaussian noise:
	\begin{equation}
		\*y = h(\*x) + \*\eta
	\end{equation}
	where $\*x$ is the state and $\*\eta \sim \mathcal{N}(0,R)$. If a prior estimate of $\*x$ is available with probability distribution $p(\*x)$, then the Bayesian posterior distribution is
	\begin{equation}
		p(\*x|\*y) \propto p(\*x)p(\*y|\*x)
		\label{eq:bayesPost}
	\end{equation}
	where $p(\*y|\*x)$ is the measurement likelihood distribution. If the measurement noise is zero-mean and Gaussian-distributed, then $p(\*y|\*x) = \mathcal{N}\big(\*y - h(\*x), R\big)$. Thus, we can rewrite Eq. \ref{eq:bayesPost} as
	\begin{equation}
		p(\*x|\*y) \propto p(\*x) \cdot \exp\Big[-\frac{1}{2}\big(\*y - h(\*x)\big)^T R^{-1} \big(\*y - h(\*x)\big)\Big].
		\label{eq:priorLik}
	\end{equation}
	
	This work addresses the approximation of $p(\*x|\*y)$ for nonlinear measurements. If the prior is Gaussian-distributed with $p(\*x) \sim \mathcal{N}(\bar{\*x}, \bar{P})$ and the measurement is \textit{linear}, then $\*x$ and $\*y$ are jointly Gaussian. The posterior distribution $p(\*x|\*y)$ is Gaussian with mean $\hat{\*x}$ and covariance $\hat{P}$ given by the usual Kalman update \cite{bar2001estimation}.  If the measurement is nonlinear, the EKF update approximates $\hat{\*x}$ and $\hat{P}$ using the measurement Jacobian $H = \frac{dh(\*x)}{d\*x}\Big|_{\*x=\bar{\*x}}$. However, the measurement Jacobian can be a poor approximation of the measurement model for highly nonlinear measurements. 
	
	To that end, we would like to gradually introduce the measurement likelihood $p(\*y|\*x)$ into Eq. \ref{eq:priorLik}, reevaluating the measurement Jacobian as we go. We begin by expressing the argument of the measurement likelihood distribution as a sum of $N$ terms:
	\begin{equation}
		p(\*x|\*y) \propto p(\*x) \cdot \exp\Big[-\frac{1}{2}\sum_{i=1}^N \frac{1}{N}\big(\*y - h(\*x))^T R^{-1} (\*y - h(\*x)\big)\Big].
		\label{eq:likSum}
	\end{equation}
	Using a convenient exponential identity, 
	\begin{equation}
		p(\*x|\*y) \propto  p(\*x) \cdot \prod_{i=1}^N \exp\Big[-\frac{1}{2}\frac{1}{N}\big(\*y - h(\*x))^T R^{-1} (\*y - h(\*x)\big)\Big].
	\end{equation}
	Finally, we move the term $N$ inside the covariance matrix inversion to form the updated expression:
	\begin{equation}
		p(\*x|\*y) \propto  p(\*x) \cdot \prod_{i=1}^N \exp\Big[-\frac{1}{2}\big(\*y - h(\*x))^T (N\!R)^{-1} (\*y - h(\*x)\big)\Big].
		\label{eq:NR}
	\end{equation}
	
	The formulation of the measurement likelihood in Eq. \ref{eq:NR} suggests that measurement information can be incorporated into the state estimate using a series of $N$ Kalman updates with an inflated measurement covariance matrix. Indeed, performing $N$ Kalman updates with measurement covariance $N\!R$ is \textit{equivalent} to performing a single Kalman update with measurement covariance $R$ for linear measurements. 
	
	We prove this equivalence in two steps. First, we show that the recursively-updated covariance $P^{(N)}$ is equivalent to the Kalman-updated covariance $\hat{P}$ for linear measurements. Then we do the same for the state update. 
	
	\begin{theorem}[Linear Covariance Convergence]
		Given prior covariance matrix $\bar{P}$ and linear measurement $\*y = H\*x$, the $N$th iterative covariance update
		\begin{equation}\label{eq:BRUF-update}
			\begin{gathered}
				P^{(i)} = (I - K^{(i)}H)P^{(i-1)},
				\quad 1 \leq i \leq N
			\end{gathered}
		\end{equation}
		is exactly equivalent to the Kalman update
		\begin{equation}
			\hat{P} = \left(I - K H\right)\bar{P} 
		\end{equation}
		where 
		\begin{equation}
			K^{(i)} = P^{(i-1)}H^T\left(H P^{(i-1)}H^T + N R\right)^{-1},\\
		\end{equation}
		and $P^{(0)} = \bar{P}$.
		\label{thm:bruf-cov}
	\end{theorem}
	\begin{proof}
		Consider the matrix inversion lemma~\cite{petersen2008matrix}:
		\begin{equation}\label{eq:woodbury}
			(P^{-1} + H^T R^{-1} H)^{-1} = \left(I - PH^T \left(HPH^T + R\right)^{-1}H\right)P.
		\end{equation}
		Using~\eqref{eq:woodbury}, it is possible to rewrite the recursive covariance update~\eqref{eq:BRUF-update} as
		\begin{equation}
			P^{(i)} = \left[\left[P^{(i-1)}\right]^{-1} + \frac{1}{N}H^T R^{-1} H\right]^{-1},\quad 1 \leq i \leq N.
			\label{eq:bruf-cov}
		\end{equation}
		Thus, by simple manipulation,
		\begin{align*}
			P^{(N)} &= \left[\left(\bar{P}\right)^{-1} + \frac{1}{N}\sum_{i=1}^N H^T R^{-1} H\right]^{-1},\\
			&= \left[\left(\bar{P}\right)^{-1} + H^T R^{-1} H\right]^{-1} = \hat{P},
		\end{align*}
		as required.
	\end{proof}
	
	The matrix inversion lemma \eqref{eq:woodbury} suggests an alternative form of the Kalman covariance update. The \textit{information matrix} is the inverse of the covariance matrix, and it is updated using the additive expression \cite{bar2001estimation}:
	\begin{equation}
		\hat{P}^{-1} = \bar{P}^{-1} + H^T R^{-1} H.
	\end{equation}
	Similarly, we can define the \textit{information state}, $\*z$, as the product of the information matrix and the state vector \cite{bar2001estimation}:
	\begin{equation}
		\*z = P^{-1}\*x.
	\end{equation}
	The information state update is: 
	\begin{equation}
		\*{\hat{z}} = \*{\bar{z}} + H^T R^{-1} \*y.
	\end{equation}
	We employ the information state representation to prove the equivalence of the BRUF state update to the usual Kalman update for linear measurements. 
	
	\begin{theorem}[BRUF Linear State Convergence] 
		Given prior covariance matrix $\bar{P}$ and linear measurement $\*y = H\*x$, the $N$th iterative state update
		\begin{equation}\label{eq:BRUF-state-update}
			\begin{gathered}
				\*x^{(i)} = \*x^{(i-1)} + K^{(i)}\left(\*y - H\*x^{(i-1)}\right),
				\quad 1 \leq i \leq N
			\end{gathered}
		\end{equation}
		is exactly equivalent to the Kalman update
		\begin{equation}
			\*{\hat{x}} = \*{\bar{x}} + K\left(\*y - H\*{\bar{x}}\right)
		\end{equation}
		where $K$ is the Kalman gain expressed in terms of the updated covariance \cite{bar2001estimation},
		\begin{equation}
			K = \hat{P}H^TR^{-1}, 
		\end{equation} 
		and $\*x^{(0)} = \*{\bar{x}}$.
		\label{thm:bruf-state-update}
	\end{theorem}
	\begin{proof}
		We begin by re-expressing \eqref{eq:BRUF-state-update} in terms of the information state \cite{bar2001estimation}:
		\begin{equation}\label{information-state-i}
			\*z^{(i)} = \left[P^{(i)}\right]^{-1}\*x^{(i)}
		\end{equation}
		where $\*z^{(0)} = \*{\bar{z}} = \bar{P}^{-1}\*{\bar{x}}$, and 
		\begin{equation}
			\*z^{(i+1)} = \*z^{(i)} + H^T(NR)^{-1}\*y
		\end{equation}
		is the information state update for recursive update step $i$.
		
		The $i$th information state update can be expressed in terms of the prior information state, $\*{\bar{z}}$, as
		\begin{equation}
			\*z^{(i+1)} = \*{\bar{z}} + \frac{i}{N}H^TR^{-1}\*y.
		\end{equation}
		Thus, 
		\begin{equation}
			\*{\hat{z}} = \*z^{(N)} = \*{\bar{z}} + H^TR^{-1}\*y.
			\label{eq:inf-state-update}
		\end{equation}
		Equation \eqref{eq:inf-state-update} is the usual expression for the information state update.
		
		Next, we can use \eqref{information-state-i} to express the updated state in terms of the updated information state:
		\begin{equation}
			\*{x}^{(N)} = P^{(N)}\*{z}^{(N)}.
		\end{equation}
		From Theorem \ref{thm:bruf-cov}, we know that
		\begin{equation}
			P^{(N)} = \hat{P} = \left[\left(\bar{P}\right)^{-1} + H^T R^{-1} H\right]^{-1}.
			\label{eq:p-plus}
		\end{equation}
		Rearranging \eqref{eq:p-plus}, it is easy to show that
		\begin{equation}
			\bar{P} = \big[\hat{P}^{-1} - H^TR^{-1}H\big]^{-1}.
		\end{equation}
		Finally, 
		\begin{equation}
			\begin{aligned}
				\*x^{(N)} &= P^{(N)}\*{z}^{(N)} \\
				&= \hat{P}\*{z}^{(N)} \\
				&= \hat{P}\*{\bar{z}} + \hat{P}H^TR^{-1}\*y \\
				&= \hat{P}(\bar{P})^{-1}\*{\bar{x}} + \hat{P}H^TR^{-1}\*y \\
				&= \hat{P}\big[(\hat{P})^{-1} - H^TR^{-1}H\big]\*{\bar{x}} + \hat{P}H^TR^{-1}\*y \\
				&= \*{\bar{x}} + \hat{P}H^TR^{-1}(\*y-H\*{\bar{x}}) \\
				&= \*{\bar{x}} + K(\*y-H\*{\bar{x}})
			\end{aligned}
		\end{equation}
		which concludes the proof.
	\end{proof}
	
	The BRUF is an approximation of the recursive linear update presented in Theorems \ref{thm:bruf-cov} and \ref{thm:bruf-state-update}. In the BRUF update, the measurement Jacobian is recomputed at the current state estimate for each update step. This allows us to compute a better representation of the Bayesian posterior than a single EKF update, which is computed in a single step using the linearization about the prior state estimate. The BRUF update is presented in Algorithm \ref{alg:bruf}.
	
	\begin{algorithm}[H]
		\begin{algorithmic}[1]
			\Require $\*{\bar{x}}$, the prior state estimate; $\bar{P}$, the prior covariance; $\*y$, a measurement; $R$, the measurement covariance; $N$, the number of steps
			\State $\*x^{(0)} \gets \*{\bar{x}}$
			\State $P^{(0)} \gets \bar{P}$
			\For{$i = 1 \dots N$} 
			\State $H^{(i)} = \frac{dh(\*x)}{d\*x}\big\rvert_{\*x = \*x^{(i-1)}}$ 
			\State $S^{(i)} = H^{(i)}P^{(i-1)}[H^{(i)}]^T + N\!R$
			\State $K^{(i)} = P^{(i-1)} [H^{(i)}]^T S^{-1}$ 
			\State $\*x^{(i)} \gets \*x^{(i-1)} + K^{(i)}\left[\*y - h(\*x^{(i-1)})\right]$ 
			\State $P^{(i)} \gets \left[I - K^{(i)}H^{(i)}\right] P^{(i-1)}$ 
			\EndFor 
			\State $\*{\hat{x}} \gets \*x^{(N)}$
			\State $\hat{P} \gets P^{(N)}$
		\end{algorithmic}
		\caption{The Bayesian Recursive Update}
		\label{alg:bruf}  
	\end{algorithm}
	
	At the beginning of each BRUF update step, the measurement Jacobian $H^{(i)}$ is computed using the current state estimate. The innovation covariance $S^{(i)}$ is computed using inflated measurement noise $N\!R$. Then, the resulting Kalman gain is used to update the state and covariance. This process is repeated $N$ times, where $N$ is an integer chosen by the designer. The BRUF update is equivalent to the EKF update for $N=1$.
	
	The BRUF is further motivated by two numerical examples. The first shows the evolution of the state and covariance matrix during the BRUF update for a Gaussian-distributed prior state estimate and a nonlinear measurement. The second demonstrates the filtering performance of the BRUF in an aerospace application.

	\subsection{Range Observation Example}
	\label{sec:range-observation-bruf}
	
	Consider prior distribution
	\begin{equation}
		p(\*x) \sim \mathcal{N}(\*{\bar{x}}, \bar{P})
	\end{equation}
	with state estimate $\*{\bar{x}} = \begin{bmatrix} -3 & 0 \end{bmatrix}^T$ and covariance matrix $\bar{P} = \begin{bmatrix} 1 & 0.5 \\ 0.5 & 1 \end{bmatrix}$. 
	
\begin{figure}[h]
		\centering
		\includegraphics[width=0.75\textwidth]{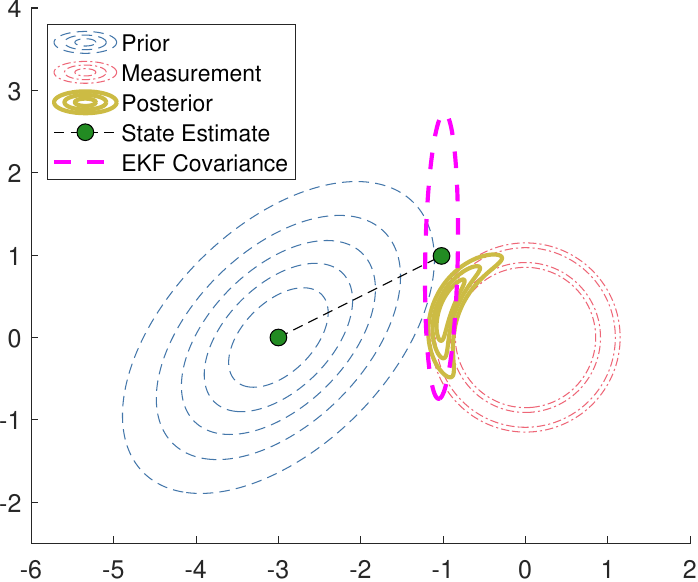}
		\caption{Result of EKF update for range observation example. A contour plot of the prior distribution is shown in blue. The prior state estimate, $\*{\bar{x}}$, lies at the mean of the prior distribution. A contour plot of the measurement likelihood distribution is shown in red. The true Bayesian posterior (yellow) was approximated using a grid method. The EKF update moves the state estimate toward the area of concentrated probability mass, but it does not represent the true statistics of the posterior distribution. The $2\sigma$ covariance ellipse of the updated distribution is shown in magenta.}
		\label{fig:banana-ekf}
	\end{figure}

	An instrument measures the distance from the origin. The measurement model is 
	\begin{equation}
		\*y = h(\*x) + \eta
	\end{equation}
	where
	\begin{equation}
		h(\*x) = \sqrt{x^2 + y^2}
	\end{equation}
	and 
	\begin{equation}
		\eta \sim \mathcal{N}(0,R)
	\end{equation}
	with $R=0.1^2$. Figure \ref{fig:banana-ekf} shows the result of the EKF update for a measurement of $\*y = 1$. Clearly, a single EKF update fails to capture the statistics of the posterior distribution. 
	
\begin{figure}[h]
		\centering
		\includegraphics[width=0.9\textwidth]{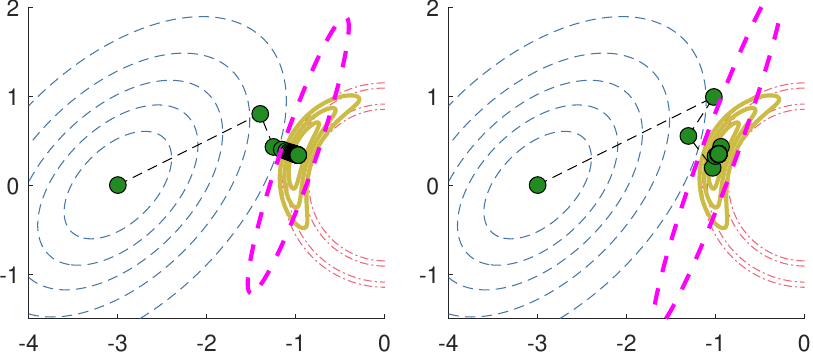}
		\caption{Results of BRUF update (left) and IEKF update (right) for range observation example. A value of $N=25$ update steps was chosen for each algorithm. The final $2\sigma$ covariance ellipses computed in each approach are shown in magenta.}
		\label{fig:banana-bruf-iekf}
	\end{figure}

	Figure \ref{fig:banana-bruf-iekf} shows the results of the BRUF update and the IEKF update for this example. The BRUF update moves the state smoothly from prior to posterior, settling near the \textit{maximum a posteriori} (MAP). Each algorithm was executed for $N=25$ update steps. The IEKF update was computed using an exact line search method \cite{havlik2015performance} with a diminishing factor of $\frac{1}{2}$. The vanilla IEKF  \cite{gelb1974applied} diverges for this problem. 
	
	Both the BRUF and the IEKF take a large first step in the direction of the EKF update. The IEKF state estimate then moves away from the MAP, before quickly settling on it. The BRUF state estimate moves smoothly toward the MAP after the second step. In fact, the larger the choice of $N$, the smoother the BRUF state update trajectory. In order to take advantage of of this smoothness without incurring the undesirable computational cost of large $N$, we introduce two modifications to the BRUF. 
	
	First, observe that Eq. \ref{eq:likSum} can be written as:
	\begin{equation}
		p(\*x|\*y) \sim C \cdot p(\*x) \cdot \exp\Big[-\frac{1}{2}\sum_{i=1}^N c_i\big(\*y - h(\*x))^T R^{-1} (\*y - h(\*x)\big)\Big]
	\end{equation}
	where 
	\begin{equation}
		\sum_{i=1}^N c_i = 1.
	\end{equation}
	In the BRUF, each coefficient $c_i$ is equal to $1/N$. However, the results in Theorems \ref{thm:bruf-cov} and \ref{thm:bruf-state-update} hold for \textit{any} set of coefficients $c_i$ as long as they sum to 1. We would like to use a set of increasing coefficients; this way, the first step (with poorly-representative measurement Jacobian $H^{(1)}$) has relatively low weight compared to later steps. The highest weights are applied at the end of the update, as the state estimate approaches the MAP. Beginning with the sum of an algebraic series of $N$ terms,
	\begin{equation}
		1 + 2 + 3 + \hdots + N = \frac{N(N+1)}{2},
		\label{eq:alg-series}
	\end{equation}
	we can write: 
	\begin{equation}
		\frac{1 + 2 + 3 + \hdots + N}{N(N+1)/2} = 1
		\label{eq:alg-series-sum-1}
	\end{equation}
	Thus, we define the \textit{Variable-Step BRUF} (VS-BRUF)\footnote{Here, we define the VS-BRUF coefficients using a sum of N terms between 1 and $N$. However, \textit{any} series with a closed-form expression for the sum could be used to compute the coefficients in a similar manner to Eqs. \ref{eq:alg-series}-\ref{eq:alg-series-sum-1}. }  by replacing the covariance $N\!R$ in each iteration of Algorithm \ref{alg:bruf} with $\frac{1}{c_i}R$, where
	\begin{equation}
		c_i = \frac{i}{N(N+1)/2}.
		\label{eq:vs-bruf}
	\end{equation}
	The VS-BRUF coefficients are similar in nature to those of the RUF \cite{zanetti2011recursive}. The VS-BRUF update is shown in Algorithm \ref{alg:vs-bruf}. Changes from the original BRUF update are highlighted in yellow.
	\begin{algorithm}[H]
		\begin{algorithmic}[1]
			\Require $\*{\bar{x}}$, the prior state estimate; $\bar{P}$, the prior covariance; $\*y$, a measurement; $R$, the measurement covariance; $N$, the number of steps; \colorbox{highlighter}{$\{c_i\}$, $i=1 \hdots N$, a set of coefficients defined according to Eq. \ref{eq:vs-bruf}}
			\State $\*x^{(0)} \gets \*{\bar{x}}$
			\State $P^{(0)} \gets \bar{P}$
			\For{$i = 1 \dots N$} 
			\State $H^{(i)} = \frac{dh(\*x)}{d\*x}\big\rvert_{\*x = \*x^{(i-1)}}$ 
			\State $S^{(i)} = H^{(i)}P^{(i-1)}[H^{(i)}]^T +$ \colorbox{highlighter}{$\frac{1}{c_i}$}$R$
			\State $K^{(i)} = P^{(i-1)} [H^{(i)}]^T S^{-1}$ 
			\State $\*x^{(i)} \gets \*x^{(i-1)} + K^{(i)}\left[\*y - h(\*x^{(i-1)})\right]$ 
			\State $P^{(i)} \gets \left[I - K^{(i)}H^{(i)}\right] P^{(i-1)}$ 
			\EndFor 
			\State $\*{\hat{x}} \gets \*x^{(N)}$
			\State $\hat{P} \gets P^{(N)}$
		\end{algorithmic}
		\caption{The Variable-Step Bayesian Recursive Update (VS-BRUF)}
		\label{alg:vs-bruf}  
	\end{algorithm}
	
	Continuing this logic, we can think of $N$ as the inverse of the length of the step that our algorithm takes: the larger the $N$, the more iterations that are performed, thus the smaller the step length is. 
	If one thinks of the update,
	\begin{equation}\label{eq:first-stage}
		x = x + \underbrace{K(\*y - h(\*x))}_{\text{First Stage}}
	\end{equation}
	as equivalent to the first stage of an explicit Runge-Kutta method, then it is possible to find an estimate of the error by utilizing an ``embedded method'' for error control. In this work, we want to make use of the minimal extra computation to compute the error, thus, for our embedded method we take the explicit midpoint rule, which is a type of predictor-corrector method. 
	
	The full method that we utilize is given by the Butcher tableau,
	\begin{equation}
		\begin{array}{r|c}\begin{matrix}0 \\ 1 \end{matrix} &
			\begin{matrix}  & \\ 1 & \phantom{0} \end{matrix} 
			\\\hline
			& \begin{matrix} 1 & 0 \\ \frac{1}{2} & \frac{1}{2}\end{matrix} 
		\end{array}
	\end{equation}
	where the top left matrix represents the relative synthetic times at which the stages are computed, the top right matrix represents the quadrature used to compute the relative stages. The first row of the bottom right matrix represents the standard update defined by~\cref{eq:first-stage}, while the second row represents the embedded method used for error control. 
	The details of the implementation are taken from known methods in time integration literature~\cite{hairerwanner}.
	The full details of the implementation can be found in~\cref{alg:ec-bruf}.
	\begin{algorithm}[H]
		\begin{algorithmic}[1]
			\Require $\*{\bar{x}}$, the prior state estimate; $\bar{P}$, the prior covariance; $\*y$, a measurement; $R$, the measurement covariance; $N$, the number of steps; $a_{tol}$ and $r_{tol}$, tolerance values; $f$, $f_{min}$, and $f_{max}$, factors
			\State $\*x^{(0)} \gets \*{\bar{x}}$
			\State $P^{(0)} \gets \bar{P}$
			\State $t_c = 0$
			\State $ds = 1/N$
			\State $i = 1$
			\\
			\While{$t_c < 1$} 
			\If{$t_c + ds > 1$}
			\State{$ds \gets 1 - t_c$}
			\EndIf
			\\
			\State $H^{(i)} = \frac{dh(\*x)}{d\*x}\big\rvert_{\*x = \*x^{(i-1)}}$ \Comment{Compute EKF update using $\*x^{(i-1)}$, $P^{(i-1)}$}
			\State $S^{(i)} = H^{(i)}P^{(i-1)}[H^{(i)}]^T + \frac{1}{ds}R$ 
			\State $K^{(i)} = P^{(i-1)} [H^{(i)}]^T [S^{(i)}]^{-1}$ 
			\State $\Delta\*x^{(i)} = K^{(i)}\left[\*y - h(\*x^{(i-1)})\right]$
			\State $\tilde{\*{x}}^{(i)} = \*{x}^{(i-1)} + \Delta\*x^{(i)}$
			\State $\tilde{P}^{(i)} = \left[I - K^{(i)}H^{(i)}\right] P^{(i-1)}$ 
			\\
			\State $H^{(i)'} = \frac{dh(\*x)}{d\*x}\big\rvert_{\*x = \tilde{\*x}^{(i)}}$ \Comment{Compute EKF update using $\tilde{\*x}^{(i)}$, $\tilde{P}^{(i)}$}
			\State $S^{(i)'} = H^{(i)'}\tilde{P}^{(i)}[H^{(i)'}]^T + \frac{1}{ds}R$
			\State $K^{(i)'} = \tilde{P}^{(i)} [H^{(i)'}]^T [S^{(i)'}]^{-1}$
			\State $\Delta\*x^{(i)'} = K^{(i)'}\left[\*y - h\left(\tilde{\*x}^{(i)}\right)\right]$
			\State $\tilde{\*{x}}^{(i)'} = \*{x}^{(i-1)} + \frac{1}{2}\left(\Delta\*x^{(i)} + \Delta\*x^{(i)'}\right)$
			\State $\*s_c = a_{tol} + \max(\text{abs}(\tilde{\*x}^{(i)}), \text{abs}(\tilde{\*x}^{(i)'}))*r_{tol}$
			\State $err = \max \left(\text{RMS}\left(\frac{1}{\*s_c}\circ\left[\tilde{\*x}^{(i)} - \tilde{\*x}^{(i)'}\right]\right)\right)$
			\\
			\If{$err > 1$}
			\State $ds \gets ds \min \left(0.9, \max \left(f_{min}, f\sqrt{\frac{1}{err}}\right)\right)$
			\State reject update and go to the top of the loop
			\EndIf
			\\
			\algstore{ec-bruf} 
		\end{algorithmic}
		\caption{The Error-Controller Bayesian Recursive Update (EC-BRUF)}
		\label{alg:ec-bruf}  
	\end{algorithm}
	\begin{algorithm}
		\begin{algorithmic}
			\algrestore{ec-bruf}
			\State $t_c = t_c + ds$
			\State $\*x^{(i)} \gets \tilde{\*x}^{(i)}$
			\State $P \gets \tilde{P}^{(i)}$
			\State $i \gets i + 1$
			\State $ds \gets ds \min \left(f_{max}, \max \left(f_{min}, f \sqrt{\frac{1}{err}}\right)\right)$
			\EndWhile 
			\\
			\State $\*{\hat{x}} \gets \*x^{(i-1)}$
			\State $\hat{P} \gets P$
		\end{algorithmic}
	\end{algorithm}

\begin{figure}[h!]
		\centering
		\includegraphics[width=0.9\textwidth]{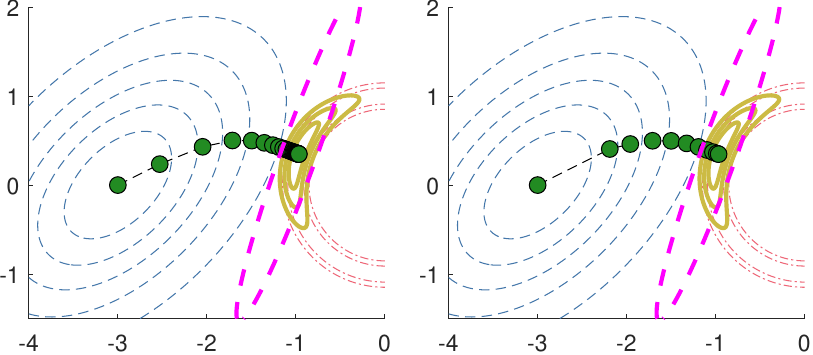}
		\caption{Results of VS-BRUF update (left) and EC-BRUF update (right) for range observation example. The VS-BRUF was executed with $N=25$. The EC-BRUF was initialized with $N=25$ but finished after 10 iterations. The final $2\sigma$ covariance ellipses computed in each approach are shown in magenta.}
		\label{fig:banana-vs-ec}
	\end{figure}

	Figure \ref{fig:banana-vs-ec} shows the results of the VS-BRUF update and the EC-BRUF update for this example. Each is very smooth compared to the BRUF and the IEKF with the same number of steps (see Fig. \ref{fig:banana-bruf-iekf}).

\begin{figure}[h]
		\centering
		\includegraphics[width=0.75\textwidth]{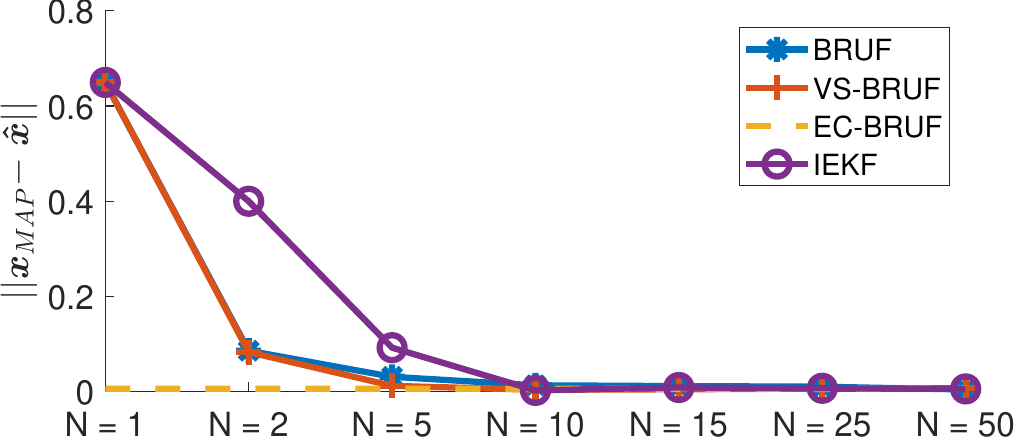}
		\caption{State estimate convergence for range observation example}
		\label{fig:banana-state-convergence}
	\end{figure}

\begin{figure}[h!]
		\centering
		\includegraphics[width=0.75\textwidth]{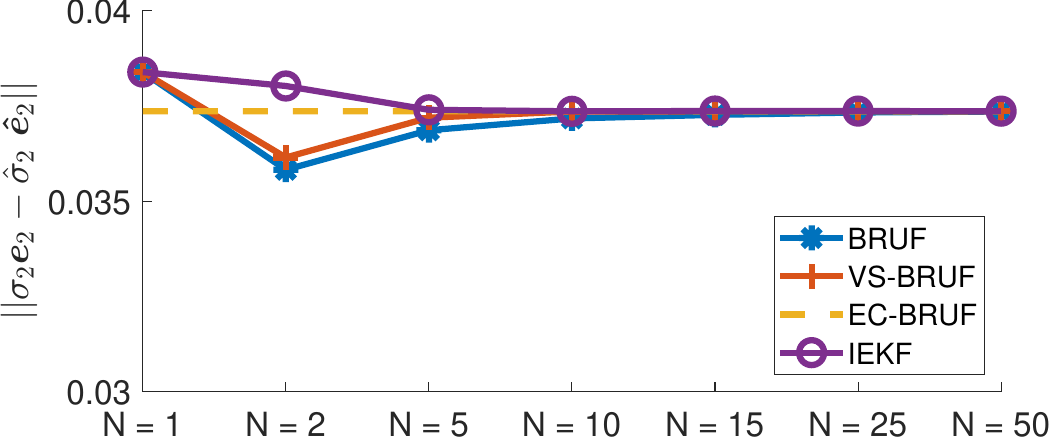}
		\caption{Covariance convergence for range observation example}
		\label{fig:banana-cov-convergence}
	\end{figure}

	A value of $N=1$ is equivalent to the EKF update for the BRUF, the VS-BRUF and the IEKF. Figure \ref{fig:banana-state-convergence} shows the convergence behavior of each approach for increasing values of $N$. The EC-BRUF was executed with $\textit{atol} = \textit{rtol} = 0.1$. It produces the same result given any value of $N$. All the algorithms reach the MAP in 10 iterations, although the BRUF variants converge faster. Also, since a line search method is necessary to run the IEKF, more computations are required per update step. Note that Figures \ref{fig:banana-state-convergence} and \ref{fig:banana-cov-convergence} represent different \textit{total} numbers of update steps $N$; they do not show error values over the course of a single measurement update. 
	
	Figure \ref{fig:banana-cov-convergence} shows the norm of the difference
	\begin{equation}
		\sigma_2 \*{e}_2 - \hat{\sigma}_2 \hat{\*{e}}_2
		\label{eq:banana-cov-plot}
	\end{equation}
	where $\sigma_2$ is the smaller eigenvalue of the true covariance of the Bayesian posterior distribution, and $\*{e}_2$ is the corresponding eigenvector. $\sigma_2$ is the length of the semi-minor axis of the covariance ellipse; the axis parallel to the range measurement. A vector difference is used to highlight the differences in \textit{angle} between the semi-minor axes; the semi-minor axis of the EKF-updated covariance ellipse in Fig. \ref{fig:banana-ekf} has a similar length to the semi-major axes of the recursively-updated covariance ellipses, but it points in the wrong direction. While there is some variation for small $N$, all the iterative approaches produce virtually the same $\hat{\sigma}_2 \hat{\*{e}}_2$ vectors for $N>2$. 
	%

	\subsection{Tracking Example}
	\label{sec:tracking-bruf}
	The so-called ``contact lens problem" in target tracking is a natural extension of the range observation problem \cite{tian2009coordinate, romeo2014particle}. The contact lens problem arises when a highly accurate range measurement is coupled with comparatively poor angular resolution. This disparity leads to a measurement likelihood distribution with probability mass concentrated near the surface of a spherical shell at distance $r$ from the receiver. This distribution mimics the shape of a contact lens. If the target is far from the receiver, the difference in uncertainty between the range direction and the plane perpendicular the range vector can be severe.
	
	Phased-array radar measurements are recorded in $r$-$u$-$v$ coordinates. The $r$-$u$-$v$ coordinate system is defined in terms of the local Cartesian coordinate system as 
	\begin{equation}
		r = \sqrt{x^2 + y^2 + z^2}
	\end{equation}
	\begin{equation}
		u = \frac{x}{r}
	\end{equation}
	\begin{equation}
		v = \frac{y}{r}
	\end{equation}
	where $u$ and $v$ are the dimensionless \textit{cosines} of the angles between the range vector and the $x$- and $y$-axes respectively. The $x$-$y$ plane lies on the face of the instrument. The measurement model is
	\begin{equation}
		\*{y}(k) = \begin{bmatrix} r(k) \\ u(k) \\ v(k) \end{bmatrix} + \*{\eta}
	\end{equation}
	where 
	\begin{equation}
		\*{\eta} \sim \mathcal{N}(\*0,R)
	\end{equation}
	and
	\begin{equation}
		R =  
		\begin{bmatrix}
			\sigma_r^2 & 0 & 0 \\
			0 & \sigma_u^2 & 0 \\
			0 & 0 & \sigma_v^2   
		\end{bmatrix}.
	\end{equation}
	
	In the simulation, a radar system tracks a target moving at nearly constant velocity \cite{bar2001estimation} far from the sensor. The state vector contains the position and velocity of the target in the local Cartesian frame.
	\begin{equation}
		\*{x} = \begin{bmatrix} x & v_x & y & v_y & z & v_z \end{bmatrix}^T
	\end{equation}
	The discrete-time dynamics propagation is
	\begin{equation}
		\*x(k+1) = F\*x(k) + \*{\nu}(k)
	\end{equation}
	where
	\begin{equation}
		F = \begin{bmatrix}
			1 & T & 0 & 0 & 0 & 0 \\
			0 & 1 & 0 & 0 & 0 & 0 \\
			0 & 0 & 1 & T & 0 & 0 \\
			0 & 0 & 0 & 1 & 0 & 0 \\
			0 & 0 & 0 & 0 & 1 & T \\
			0 & 0 & 0 & 0 & 0 & 1
		\end{bmatrix},
	\end{equation}
	and T is the length of time between measurements. The process noise is 
	\begin{equation}
		\*{\nu(k)} \sim \mathcal{N}(\*{0}, Q)
	\end{equation}
	where
	\begin{equation}
		Q = \begin{bmatrix}
			T^3/3 & T^2/2 & 0 & 0 & 0 & 0 \\
			T^2/2 & T & 0 & 0 & 0 & 0 \\
			0 & 0 & T^3/3 & T^2/2 & 0 & 0 \\
			0 & 0 & T^2/2 & T & 0 & 0 \\
			0 & 0 & 0 & 0 & T^3/3 & T^2/2 \\
			0 & 0 & 0 & 0 & T^2/2 & T 
		\end{bmatrix} \tilde{q},
	\end{equation}
	and $\tilde{q} = 10^{-4}$ m\textsuperscript{2}/s\textsuperscript{3}.
	
	The initial state of the target is 
	\begin{equation}
		\*x(0) = 
		\begin{bmatrix} 1100\;\text{km} \\
			-2\;\text{km/s} \\
			1100\;\text{km} \\ 
			-2\;\text{km/s} \\
			1100\;\text{km} \\
			-1\;\text{km/s} 
		\end{bmatrix}
	\end{equation}
	in the Cartesian frame. Radar measurements are received every $T = 1$ second. The measurement noise characteristics are $\sigma_r = 2.5$ m and $\sigma_u = \sigma_v = 1$ msin. The filter is initialized using the first two measurements. The position is estimated using:
	\begin{align}
		\hat{x}(k) = \tilde{u}(k) \tilde{r}(k) \\
		\hat{y}(k) = \tilde{v}(k) \tilde{r}(k) \\
		\hat{z}(k) = \tilde{r}(k) \sqrt{1 - \tilde{u}(k)^2 - \tilde{v}(k)^2}
	\end{align}
	where $\*y(k) = \begin{bmatrix} \tilde{r}(k) & \tilde{u}(k) & \tilde{v}(k) \end{bmatrix}^T$ is a noisy measurement, and $k = \{1,2\}$. The covariance of each measurement is placed in the Cartesian frame using the first-order approximation in \cite{tian2009coordinate}. The estimated velocity at time step 2 is \cite{bar2001estimation}:
	\begin{align}
		\hat{v}_x(2) = \frac{\hat{x}(2) - \hat{x}(1)}{T} \\
		\hat{v}_y(2) = \frac{\hat{y}(2) - \hat{y}(1)}{T} \\
		\hat{v}_z(2) = \frac{\hat{z}(2) - \hat{z}(1)}{T}
	\end{align}
	Since the first two measurements are independent, the off-diagonal terms in the initial covariance matrix $\hat{P}(2)$ are straightforward to compute using the converted position measurement covariances. The first measurement update occurs at the third time step. The target is tracked for 300 seconds.
	
\begin{figure}[h]
		\centering
		\includegraphics[width=0.9\textwidth]{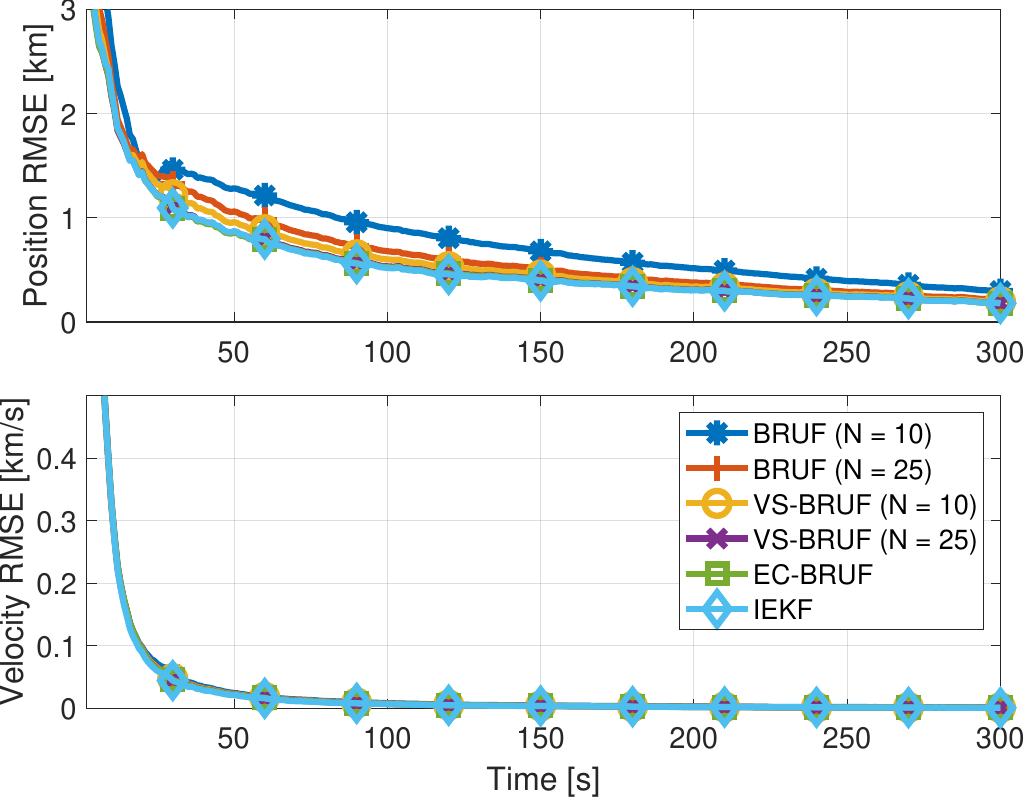}
		\caption{Tracking performance of the BRUF, the VS-BRUF, the EC-BRUF, and the IEKF. All filters quickly converge to the correct velocity. The BRUF can track the target with $N=10$ and improves for higher values of $N$. The tracking performance of the IEKF, the VS-BRUF with $N=25$, and the EC-BRUF are virtually indistinguishable.}
		\label{fig:tracking-rmse}
	\end{figure}

	Figure \ref{fig:tracking-rmse} shows the mean tracking performance of the BRUF and the IEKF over 100 Monte Carlo runs. The IEKF does not require a line search to solve this problem; iterations stopped when the difference between consecutive updates was less than $10^{-9}$, or at a maximum of 25 iterations. The error-controller adaptive BRUF was initialized with $N = 25$ and executed with $atol = rtol = 10^{-7}$. 
	
\begin{figure}[h]
		\centering
		\includegraphics[width=0.72\textwidth]{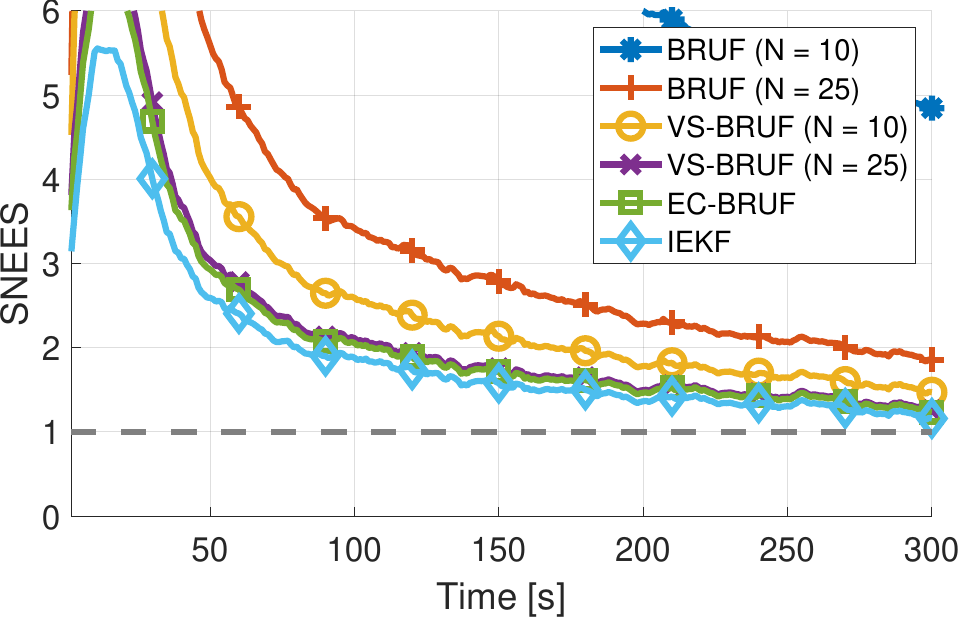}
		\caption{SNEES values for the tracking filters. In general, consistency improves with more iterations. While the BRUF produces acceptable state estimates, it is less consistent than the other filters. The IEKF, the VS-BRUF with $N=25$, and the EC-BRUF all approach a SNEES value of 1 after convergence.}
		\label{fig:tracking-snees}
	\end{figure}

	The consistency of each filter can be evaluated using the state-normalized estimation error squared (SNEES):
	\begin{equation}
		\text{SNEES}(k) = \frac{1}{n_m} \sum_{i=1}^{n_m} \frac{1}{n} \*\epsilon(k)^T P(k)^{-1} \*\epsilon(k)
	\end{equation}
	where $n$ is the length of the state vector, $n_m$ is the number of Monte Carlo runs, and $\*\epsilon(k)$ is the state estimation error at time step $k$. A SNEES value of 1 indicates that the filter is consistent; the covariance matrix is representative of the uncertainty in the state estimate.  Figure \ref{fig:tracking-snees} shows the evolution of the SNEES values for the filters in Figure \ref{fig:tracking-rmse}. 
	
	Finally, we present a runtime comparison for the filters in this example. Table \ref{tab:tracking-runtime} shows the total time it took to run 100 Monte Carlo simulations for each filter on a standard laptop. Each filter processes measurements much faster than they would arrive in a real-world tracking scenario. The processing time is contrasted with the time-averaged position RMSE:
	\begin{equation}
		\text{RMSE} = \frac{1}{n_t} \sum_{k=1}^{n_t} \sqrt{\frac{1}{n_m}\sum_{i=1}^{n_m} ||\*\epsilon_p(k)||_2^2}
	\end{equation}
	where $n_t$ is the number of time steps and $\*\epsilon_p(k)$ is the vector of position errors at time step $k$.
	\begin{table}[h]
		\centering
		\begin{tabular}{|c|c|c|}
			\hline
			Filter & Runtime [s] & Time-Avg. Posn. RMSE [km] \\
			\hline
			EKF & 0.5 & - \\
			BRUF (N = 10) & 2.9 & 0.87  \\
			BRUF (N = 25) & 6.9 & 0.71 \\
			VS-BRUF (N = 10) & 3.0 & 0.65 \\
			VS-BRUF (N = 25) & 7.1 & 0.60 \\
			EC-BRUF (tol = $10^{-7}$) & 80.4 & 0.59 \\
			IEKF & 2.1 & 0.59 \\
			\hline
		\end{tabular}
		\caption{Runtime for 100 Monte Carlo trials of each tracking filter. Each run contains 300 measurement updates. Although the EKF is prone to divergence, its runtime is shown as a baseline for comparison. The IEKF is the fastest, followed closely by the BRUF (N = 10) and the VS-BRUF (N = 10). The best positioning performance is given by the VS-BRUF (N = 25), the EC-BRUF, and the IEKF. The EC-BRUF has a long runtime for this problem, since tight tolerance values are required to match the accuracy of the IEKF and the VS-BRUF (N = 25).}
		\label{tab:tracking-runtime}
	\end{table}
	
	The following section presents the Bayesian Recursive Update Ensemble Kalman Filter (BRUEnKF), a modification of the Bayesian Recursive Update for ensemble filters. We draw comparisons between the BRUEnKF and existing \textit{particle flow filters}, which smoothly move states from the prior distribution to the Bayesian posterior. These filters can solve highly nonlinear problems where single-state filters fail. Finally, we demonstrate the performance of the BRUEnKF on a high-dimensional problem. 
	
	\section{The Bayesian Recursive Update EnKF}
	\label{sec:bruenkf}
	
	The BRUEnKF update begins with an \textit{ensemble} of state estimates that represent the prior distribution. The BRUEnKF moves the ensemble states through the BRUF update together, recomputing the ensemble covariance as it evolves. 
	
	In general, the EnKF \cite{evensen2009data} is useful for high-dimensional systems when computing and inverting large state covariance matrices at each time step is impractical. 
	%
	%
	In the \textit{linearized EnKF}, a Kalman gain is computed for each ensemble member using:
	\begin{equation}\label{eq:linearized-EnKF}
		K_j = \Cov(\bar{\*X},\bar{\*X})H_j^T\left(H_j \Cov(\bar{\*X},\bar{\*X}) H_j^T + R \right)^{-1}
	\end{equation}
	where $\bar{\*{X}}$ is the ensemble of prior state estimates $\bar{\*x}_j$, $j = 1 \hdots M$, and $H_j = \frac{dh(\*x)}{d\*x}\big|_{\*x=\bar{\*x}_j}$ is the measurement Jacobian evaluated at state $\bar{\*x}_j$. Both the statistical EnKF and the linearized EnKF perform the measurement update in one step. Algorithm \ref{alg:enkf} outlines the linearized EnKF update. The EnKF update begins by \textit{inflating} the ensemble members by factor $\alpha$. Then, an EKF update is performed on each ensemble member using a predicted measurement perturbed by an artificial noise. These measures provide statistical convergence guarantees \cite{popov2020explicit}. The state estimate at any time step is simply the mean of the updated ensemble states.
	\begin{algorithm}[H]
		\begin{algorithmic}[1]
			\Require $\bar{\*X}$, an ensemble of $M$ samples $\bar{\*x}_j$ from the prior distribution, where $j=1 \hdots M$; $\alpha$, an inflation factor
			\State $\*m = \frac{1}{M}\sum_{j=1}^M \*x_j$ \Comment{Recompute mean}
			\State $\*X \gets \*m + \alpha(\bar{\*X} - \*m)$ \Comment{Perform inflation}
			\State $P = \frac{1}{M-1}(\bar{\*X} - \*m)(\bar{\*X} - \*m)^T$ \Comment{Recompute covariance}
			\For{$j = 1 \hdots M$} 
			\State $\*{\hat{y}}_j = h(\*x_{j}) + \gamma_{j}$ \Comment{$\gamma_j \sim \mathcal{N}(\*0, R)$}
			\State $H_j = \frac{\partial h}{\partial\*x}\big\rvert_{\*x = \*x_{j}}$
			\State $S = H_j P H_j^T + R$
			\State $K = P H_j^T S^{-1}$
			\State $\hat{\*x}_{j} \gets \bar{\*x}_{j} + K(\*y - \*{\hat{y}}_j)$ \Comment{Build $\hat{\*X}$}
			\EndFor
		\end{algorithmic}
		\caption{The Linearized EnKF Update}
		\label{alg:enkf}  
	\end{algorithm}
	
	Algorithm \ref{alg:bruenkf} presents the BRUEnKF. The approach is an extension of the linearized EnKF. In the BRUEnKF, a BRUF update is performed on each ensemble member using ensemble covariance $P^{(i-1)}$. The inflation factor $\alpha^{1/N}$ is applied and the ensemble covariance is recomputed at each BRUF update step. 
	\begin{algorithm}[H]
		\begin{algorithmic}[1]
			\Require $\bar{\*X}$, an ensemble of $M$ samples $\bar{\*x}_j$ from the prior distribution; $N$, the number of BRUF steps; $\alpha$, an inflation factor
			\State $\*X^{(0)} \gets \bar{\*X}$
			\For{$i = 1 \hdots N$} 
			\State $\*m^{(i-1)} = \frac{1}{M}\sum_{j=1}^M \*x_j^{(i-1)}$ 
			\State $\*X^{(i-1)} \gets \*m^{(i-1)} + \alpha^{1/N}\left[\*X^{(i-1)} - \*m^{(i-1)}\right]$ 
			\State $P^{(i-1)} = \frac{1}{M-1}\left[\*X^{(i-1)} - \*m^{(i-1)}\right]\left[\*X^{(i-1)} - \*m^{(i-1)}\right]^T$ 
			\For{$j = 1 \hdots M$} 
			\State $\*{\hat{y}}_j = h\left(\*x^{(i-1)}_j\right) + \gamma_{j}$ \Comment{$\gamma_j \sim \mathcal{N}(\*0, R)$}
			\State $H_j = \frac{\partial h}{\partial\*x}\big\rvert_{\*x = \*x^{(i-1)}_j}$
			\State $S = H_j P^{(i-1)} H_j^T + NR$
			\State $K = P^{(i-1)} H_j^T S^{-1}$
			\State $\*x^{(i)}_j \gets \*x^{(i-1)}_j + K(\*y - \*{\hat{y}}_j)$ \Comment{Build $\*X^{(i)}$}
			\EndFor
			\EndFor
			\State $\hat{\*X} \gets \*X^{(N)}$
		\end{algorithmic}
		\caption{The Bayesian Recursive Update EnKF (BRUEnKF) Update}
		\label{alg:bruenkf}  
	\end{algorithm}
	
\begin{figure}[h]
		\centering
		\includegraphics[width=0.75\textwidth]{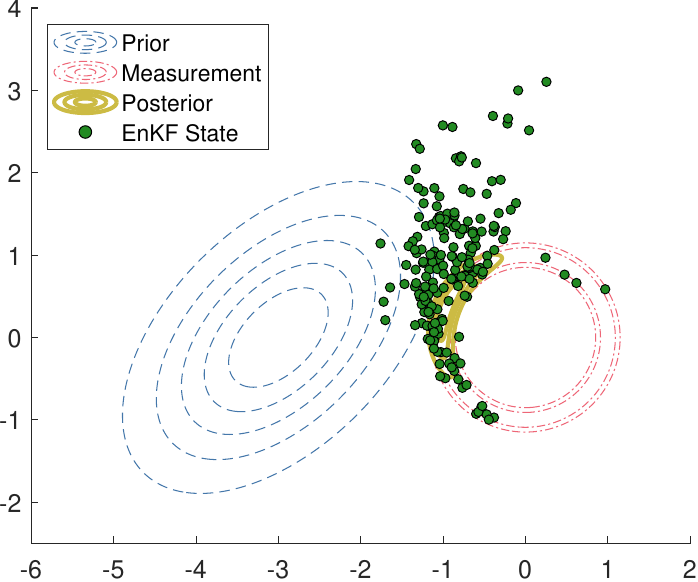}
		\caption{Results of linearized EnKF update for range observation example with 200 ensemble states (green). Before the update, the ensemble is distributed according to the prior distribution (blue). The linearized EnKF update (Algorithm \ref{alg:enkf}) is applied to each ensemble member. While some ensemble members reach the measurement likelihood distribution (red), most end up far from the MAP.}
		\label{fig:banana-enkf}
	\end{figure}

	To give a more intuitive visual understanding of the BRUEnKF update, we revisit the range observation example from Section \ref{sec:range-observation-bruf}. Fig. \ref{fig:banana-enkf} shows the result of the linearized EnKF update. Like the EKF update, the linearized EnKF update fails to capture the statistics of the Bayesian posterior.
	
\begin{figure}
		\centering
		\includegraphics[width=\textwidth]{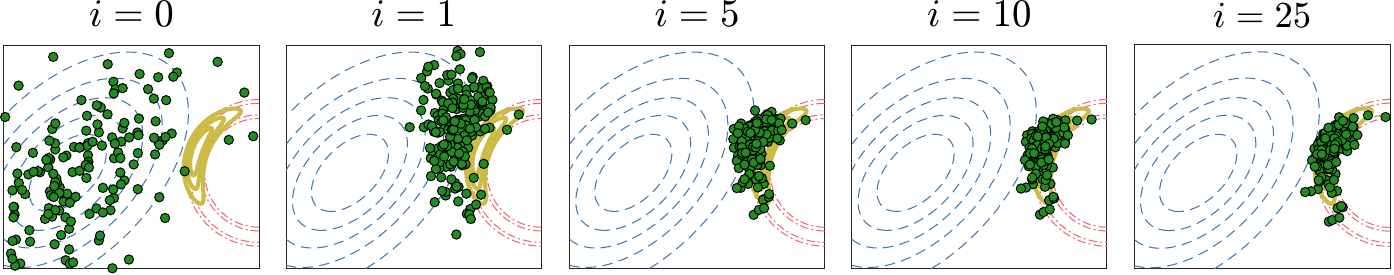}
		\caption{Results of BRUEnKF update for range observation example. The ensemble begins distributed according to the Gaussian prior distribution (far left). Much like the BRUF update with respect to the EKF, the ensemble resembles the EnKF update after the first iteration. The change in the distribution of the ensemble members is much greater between $i=1$ (center left) and $i=5$ (center) than between $i=10$ (center right) and $i=25$ (right).}
		\label{fig:banana-bruenkf}
	\end{figure}

	Figure \ref{fig:banana-bruenkf} shows the evolution of the BRUEnKF update for $N = 25$ update steps. The ensemble settles nicely on the crescent-moon shape of the posterior distribution. However, like the BRUF update, the motion of the ensemble members is faster at the beginning of the update than the end. We introduce the VS-BRUEnKF and the EC-BRUEnKF to smooth out the BRUEnKF update. The VS-BRUEnKF is identical to the BRUEnKF, except that $\alpha^{c_i}$ replaces $\alpha^{1/N}$ on line 4 of Algorithm \ref{alg:bruenkf}, and $1/c_i$ replaces $N$ on line 9. The coefficients $c_i$ are generated using Equation \ref{eq:vs-bruf}. In a similar vein, the EC-BRUEnKF extends Algorithm \ref{alg:ec-bruf} for an ensemble update. A full description of the implementation is omitted for the sake of brevity.
	
	%

\begin{figure}
		\centering
		\includegraphics[width=\textwidth]{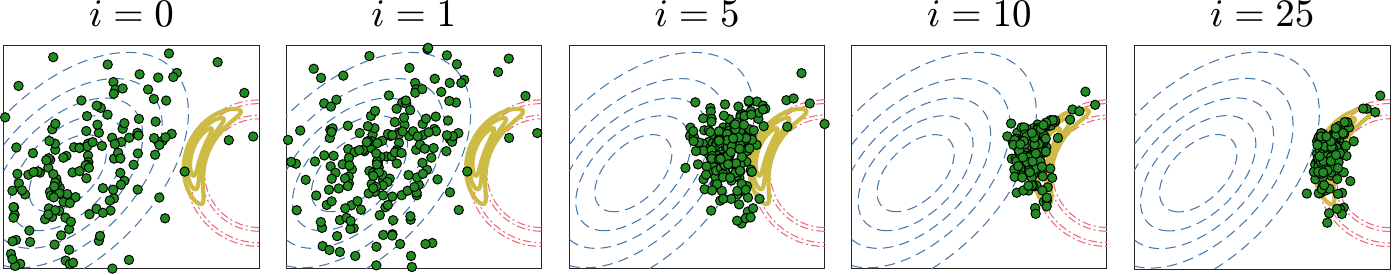}
		\caption{Results of VS-BRUEnKF update for range observation example}
		\label{fig:banana-vs-bruenkf}
	\end{figure}

\begin{figure}
		\centering
		\includegraphics[width=\textwidth]{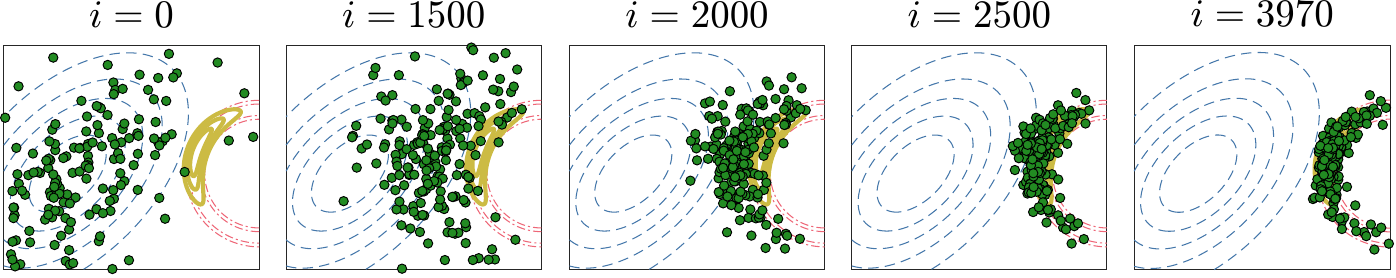}
		\caption{Results of EC-BRUEnKF update for range observation example}
		\label{fig:banana-ec-bruenkf}
	\end{figure}

	Figure \ref{fig:banana-vs-bruenkf} shows the evolution of the VS-BRUEnKF for $N=25$ update steps. The evolution of the ensemble is slowed significantly compared to the BRUEnKF update. The first update step induces very little motion. By the fifth step, the ensemble has collected into a shape resembling the posterior distribution. At the end of the update, the distribution tightens and moves onto the posterior. The EC-BRUEnKF shows similar behavior in Figure \ref{fig:banana-ec-bruenkf}. The end result is an excellent representation of the Bayesian posterior distribution. The EC-BRUEnKF runs 3970 iterations for this problem with $a_{tol} = r_{tol} = 10^{-6}$. 
	
	We demonstrate the performance of these ensemble filters on a Lorenz '96 system with a highly nonlinear measurement.
	
	\section{Lorenz '96}
	\label{sec:bruenkf-l96}
	
	The Lorenz '96 dynamics \cite{lorenz1996predictability, lorenz1998optimal} are defined by
	\begin{equation}
		\dot{x}_i = (x_{i+1} - x_{i-2})x_{i-1} - x_i + F
		\label{eq:l96}
	\end{equation}
	for $i = 1, \hdots, n$, where $n$ is the dimension of the state vector. We make the standard choices of $n = 40$ and $F = 8$ for a chaotic dynamical system. The state vector represents a 1-dimensional circular domain; hence, the dynamics of the boundary terms are defined using $x_0 = x_{40}$, $x_{-1} = x_{39}$, and $x_{41} = x_1$ \cite{asch2016data}.
	
	The measurement model is $\*y = h(x)$:
	\begin{equation}
		h(x) = \frac{x}{2}\left[1 + \left(\frac{|x|}{f}\right)^{\gamma-1}\right]
		\label{eq:l96-meas}
	\end{equation}
	where $f=10$ scales the measurement, and $\gamma$ controls the nonlinearity of the measurement function, with the choice of $\gamma=1$ giving a linear measurement. We choose $\gamma = 5$ for a moderately nonlinear measurement. The measurement vector includes a nonlinear measurement of every other state value:
	\begin{equation}
		\*y = \begin{bmatrix} 
			h(x_2) \\ 
			h(x_4) \\
			\vdots \\
			h(x_{40})
		\end{bmatrix} + \*\eta
	\end{equation}
	where 
	\begin{equation}
		\*\eta \sim \mathcal{N}\left(\*0, R \right)
	\end{equation}
	is a Gaussian-distributed measurement noise vector with $R = I_{20 \times 20}$.
	
\begin{figure}
		\centering
		\includegraphics[width=0.75\textwidth]{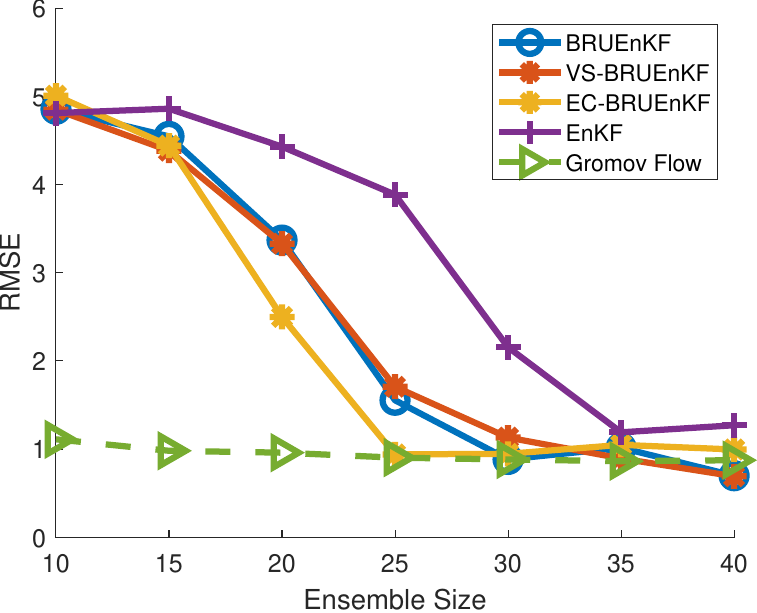}
		\caption{Time-averaged RMSE results for Lorenz '96 system with nonlinear measurement. Each point on the plot shows the mean of 10 Monte Carlo runs.}
		\label{fig:l96-rmse}
	\end{figure}

	We compare the performance of the BRUEnKF, the VS-BRUEnKF, and the EC-BRUEnKF to that of the linearized EnKF (Algorithm \ref{alg:enkf}) and Gromov flow. Gromov flow is a particle flow filter that uses a companion EKF or unscented Kalman filter (UKF) to compute the Kalman gain during the update steps and resample after the update is completed \cite{daum2016gromov}. We choose the EKF for this implementation, since the UKF requires a large number of sigma points for high-dimensional states. For a state vector of length 40, the UKF requires 80 or 81 sigma points! In contrast, we present estimation error results for ensemble sizes between 10 and 40 in Figure \ref{fig:l96-rmse}. 
	
	Gromov flow was derived as an approximation to the solution of a stochastic differential equation defined by a homotopy between the prior distribution and the Bayesian posterior \cite{daum2010exact, daum2016gromov, crouse2019consideration}. A single step of the Gromov flow update is \cite{crouse2019consideration}:
	\begin{equation}
		\*x^{(\lambda + \Delta)}_j = \*x^{(\lambda)}_j + \*f(\*x^{(\lambda)}_j,\lambda) \Delta + B \tilde{\*w}
		\label{eq:gromov-flow}
	\end{equation}
	where $\*x^{(0)}_j$ lies on the prior distribution and $\*x^{(1)}_j$ lies on the Bayesian posterior.  For $N$ update steps, the step size is $\Delta = 1/N$. The function $\*f(\*x^{(\lambda)}_j,\lambda)$ is 
	\begin{equation}
		\*f = - \left(P^{-1} + \lambda H^T R^{-1} H\right)^{-1}H^T R^{-1} \left(h\left(\*x^{(\lambda)}_j\right) - \*y \right)
	\end{equation}
	where $P$ is the EKF covariance, and $H$ is the measurement Jacobian computed at $\*x^{(\lambda)}_j$. The diffusion term $\tilde{\*w}$ is zero-mean with covariance $\Delta I$. The matrix $B$ is defined using $Q^{(\lambda)} = BB^T$, where:
	\begin{equation}
		Q^{(\lambda)} = \left(P^{-1} + \lambda H^T R^{-1} H \right)^{-1} H^T R^{-1} H \left(P^{-1} + \lambda H^T R^{-1} H \right)^{-1}.
	\end{equation}
	Intriguingly, although they were derived differently, the Gromov flow update (\ref{eq:gromov-flow}) is similar to the BRUEnKF update. The term $\*f(\*x^{(\lambda)},\lambda) \Delta$ is the same as the BRUF state update, where the Kalman gain is defined in terms of the updated covariance. Instead of the diffusion term $B \tilde{\*w}$, the BRUEnKF uses the artificial measurement noise $\gamma_j$ to ensure sufficient spread in the updated states. The BRUEnKF also inflates the ensemble and recomputes the ensemble covariance at each step, whereas Gromov flow uses the prior EKF covariance throughout. Details on the implementation of particle flow filters are provided in \cite{crouse2019consideration, ding2012implementation, pal2019particle}.
	
	Figure \ref{fig:l96-rmse} shows the time-averaged RMSE results for the five filters on a Lorenz '96 system with increments of 0.05 time steps between measurements and a total of 350 time steps. The ensembles were initialized using a distribution centered on the true initial state with covariance $P_0 = I_{40\times 40}$. Errors were calculated after a burn-in period of 50 time steps. 
	
	The dynamics were integrated using fourth-order Runge-Kutta. The inflation factor was $\alpha = 1.06$. The BRUEnKF and the VS-BRUEnKF were executed with $N = 25$ update steps. The EC-BRUEnKF parameters were with $a_{tol} = r_{tol} = 10^{-3}$, $f = 0.38^{1/2}$, $f_{max} = 6$, and $f_{min} = 0.2$. The Gromov flow filter was executed with a companion EKF. The EKF covariance was propagated using the Jacobian of the state dynamics \eqref{eq:l96}. An artificial process noise covariance $Q = 0.1\,I_{40\times 40}$ was added to the propagated EKF covariance matrix at each time step for filter health. No artificial process noise was added to any other filters. In the Gromov flow filter, resampling was conducted after each measurement update. In the resampling process, an entirely new set of particles is drawn from the distribution $\mathcal{N}(\*m,\hat{P}_{EKF})$, where $\*m$ is the mean of the updated particles and $\hat{P}_{EKF}$ is the updated EKF covariance. 
	
	Impressively, Gromov flow is able to track the system with only $M = 10$ particles. However, it is important to emphasize that the companion EKF carries a $40\times 40$ full-rank covariance matrix. This does not extend to large-scale systems; carrying a full-rank covariance matrix is essentially equivalent to carrying $n_x + 1$ ensemble members, where $n_x$ is the dimension of the state vector. Thus, Gromov flow is treated as a baseline for this example. The EC-BRUEnKF converges with $M = 25$ ensemble states. The BRUEnKF and the VS-BRUEnKF converge with $M = 30$ ensemble states. The EnKF converges with $M = 35$ ensemble states. 
	
\begin{figure}
		\centering
		\includegraphics[width=0.75\textwidth]{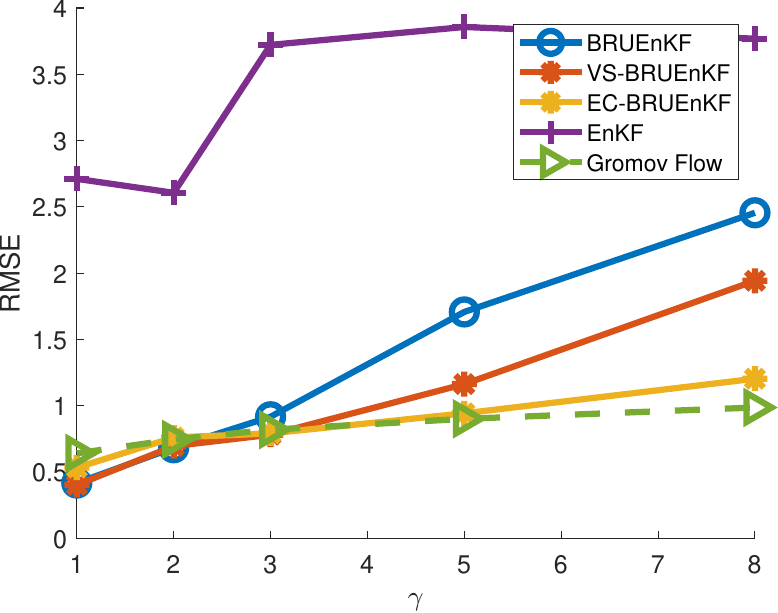}
		\caption{Time-averaged RMSE results for Lorenz '96 system with measurement nonlinearity defined by increasing values of $\gamma$ (see Equation \ref{eq:l96-meas}). Each point on the plot is the mean of 10 Monte Carlo runs.}
		\label{fig:l96-gamma}
	\end{figure}

	Figure \ref{fig:l96-gamma} shows the time-averaged RMSE results for different values of $\gamma$. Each filter carries $M=25$ states. The EnKF struggles with this system, even for linear measurements. As expected, the estimation error rises for the proposed filters and Gromov flow as $\gamma$ increases. The EC-BRUEnKF approaches the performance of Gromov flow for highly nonlinear measurements, followed by the VS-BRUEnKF and the EC-BRUEnKF.
	
	\section{Quasi-Geostrophic Model}
	\label{sec:bruenkf-qg}
	The Quasi-Geostrophic equations~\cite{foster2013finite,ferguson2008numerical,MW06,greatbatch2000four} model the flow of the streamfunction $\psi$ of some fluid on the surface of the Earth. Specifically, in terms of the vorticity $\omega$, the quadratic partial differential equation of the flow is,
	\begin{equation}\label{eq:QG}
		\begin{gathered}
			\omega_t + J(\psi,\omega) - {Ro}^{-1}\, \psi_x = {Re}^{-1}\, \Delta\omega + {Ro}^{-1}\,F, \\
			J(\psi,\omega) \equiv \psi_y\, \omega_x - \psi_x\, \omega_y,\quad \omega = -\Delta\psi,
		\end{gathered}
	\end{equation}
	where $Re=450$ is the Reynolds number, $Ro=0.0036$ is the Rossby number, $J$ is the quadratic nonlinear Jacobian term, and $F = \sin\left(\pi(y-1)\right)$ is the symmetric double gyre forcing term. In this formulation which was created for~\cite{Popov2021a,Popov2021b,popov2022model}, the spatial domain is $[0,1]\times[0,2]$ with homogeneous Dirichlet boundary conditions, and discretized with $32$ points in the $x$ direction and $64$ points in the $y$ direction. There are 2048 state values.
	The implementation used in this work is from the ODE Test Problems suite~\cite{otp,otpsoft}.
	
	For the state estimation problem, the time between observations is taken to be $0.0109$ time units, which is equivalent to about one day for the particular values of the parameters. A sparse observation of 50 evenly spaced points is utilized with error covariance $R = I_{50\times50}$. The inflation factor is $\alpha = 1.1$. The nonlinear observation operator is the same as in the previous example (Equation \ref{eq:l96-meas}), again with $\gamma = 5$. Each filter runs for 250 time steps. Errors are computed after a spinup of 50 time steps. A value of $N = 25$ update steps was chosen for the BRUEnKF and the VS-BRUEnKF. The EC-BRUEnKF is not included in this example, since it is impractical for use in large state spaces. 
	
\begin{figure}
		\centering
		\includegraphics[width=0.75\textwidth]{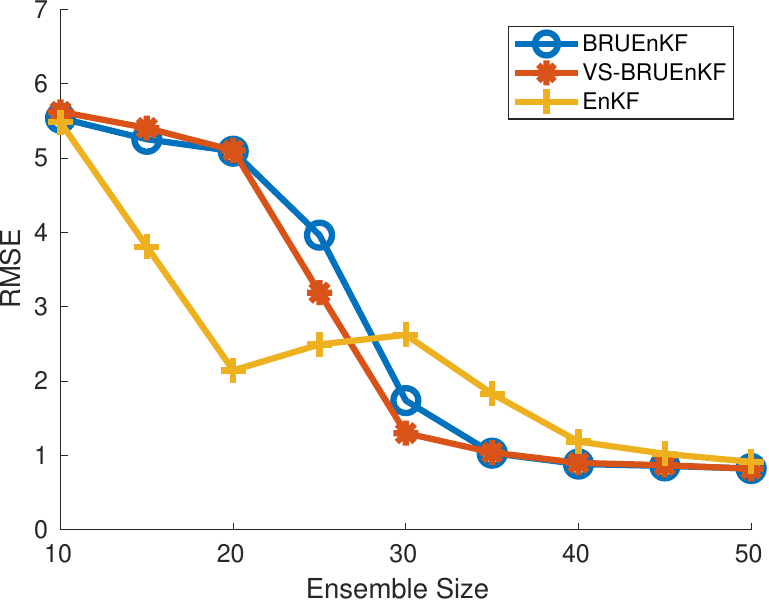}
		\caption{Time-averaged RMSE results for QG system with nonlinear measurement model defined by Equation \ref{eq:l96-meas}. Each point on the plot is the mean of 10 Monte Carlo runs.}
		\label{fig:rmse-qg}
	\end{figure}

\begin{figure}
		\centering
		\includegraphics[width=\textwidth]{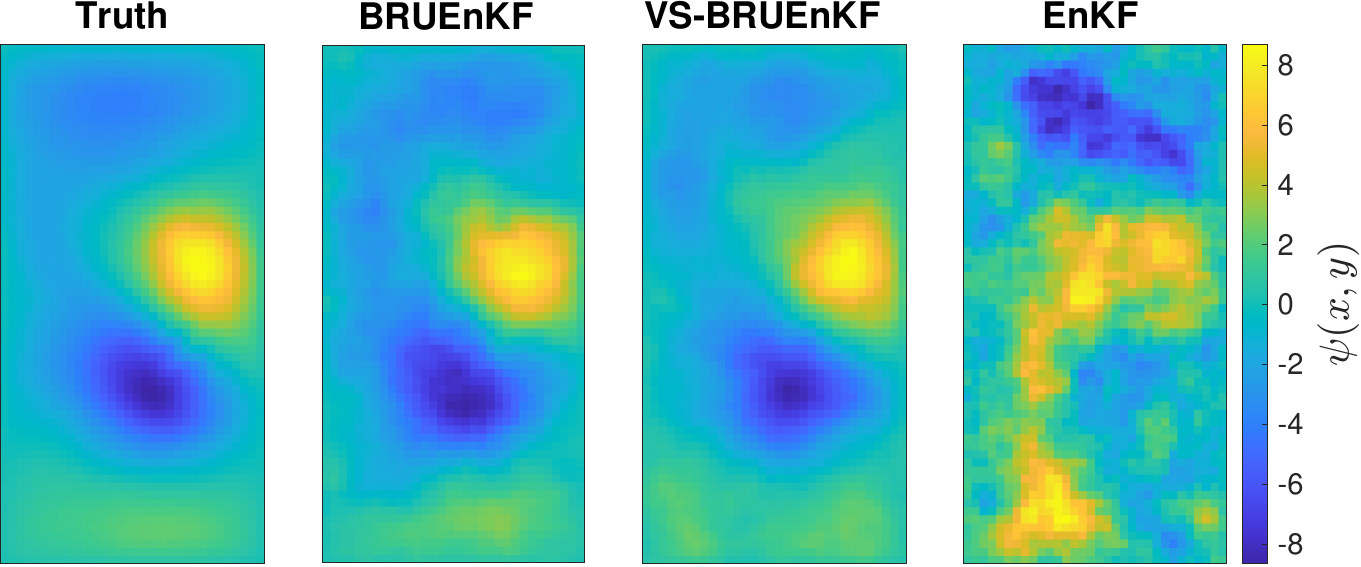}
		\caption{State estimate at final time step for BRUEnKF (center left), VS-BRUEnKF (center right), and EnKF (right) for ensemble size $M = 40$. The EnKF deviates significantly from the true state.}
		\label{fig:error-qg}
	\end{figure}

	Figure \ref{fig:rmse-qg} shows the RMSE results of the BRUEnKF, the VS-BRUEnKF, and the EnKF . The BRUEnKF and the VS-BRUEnKF show similar convergence behavior; both converge with 35--40 ensemble members. The EnKF requires 45--50 ensemble members to converge. Figure \ref{fig:error-qg} shows an example of each filter's state estimate at the final time step for an ensemble size of 40. The BRUEnKF and the VS-BRUEnKF estimates are both good representations of the true state. The EnKF errors are large. 
	
	\section{Conclusion}
	\label{sec:conclusion}
	This article presents two recursive update schemes for nonlinear measurements: the BRUF and the BRUEnKF. The recursion serves to relax the measurement linearity assumptions of the EKF and the EnKF respectively. By dividing the measurement update into $N$ steps, the update smoothly moves the state (or ensemble of states) from the prior to the posterior distribution. This work also proposes two new approaches to the Bayesian recursive update. In the first, the VS-BRUF, earlier update steps are down-weighted by over-inflating the measurement covariance; later steps are trusted more as the state vector settles on the posterior distribution. This concept is further extended in a second approach, the EC-BRUF, which uses an error controller to determine the step size.
	
	The BRUF takes its place among a class of filters that includes its predecessor, the RUF \cite{zanetti2011recursive}, and the IEKF \cite{gelb1974applied}. Unlike the IEKF, it does not iterate until a convergence criterion is met; rather, it executes a user-defined number of EKF updates for every measurement. A range observation example provides compelling visual evidence that, like the IEKF, the BRUF state estimate converges to the MAP (see Section \ref{sec:range-observation-bruf}). A convergence proof is left as future work. The introduction of the VS-BRUF and the EC-BRUF also opens the door to future avenues; any set of coefficients $1/c_i$ where $\sum_{i=1}^N c_i = 1$ can be applied to the measurement covariance matrix and still satisfy the requirements of the linear convergence proofs presented in Section \ref{sec:bruf}. 
	
	The BRUEnKF, the VS-BRUEnKF, and the EC-BRUEnKF generalize the algorithms presented in \ref{sec:bruf} for ensemble filters. Ensemble filters, like particle filters, carry a number of state vectors in order to maintain a more accurate representation of the uncertainty in the state estimate than a Gaussian distribution. Ensemble filters have the added benefit of low-order representation; instead of carrying an $n_x \times n_x$ covariance matrix, ideally, only a relatively small number of state vectors is required. The proposed ensemble filters are shown to reduce the ensemble size required for convergence in two large-scale systems (see Sections \ref{sec:bruenkf-l96} and \ref{sec:bruenkf-qg}). Future work will explore extensions of the BRUF update for other nonlinear filters; particularly particle flow filters.

	\section*{Acknowledgments}
	This work was sponsored in part by the Air Force Office of Scientific Research (AFOSR) under award number: FA9550-22-1-0419.

	\bibliographystyle{ieeetr}
	\bibliography{bib/bibliography, bib/timeintegration, bib/andreybib}

	
	
	
	
	
	
	
\end{document}